\documentclass{article}
\usepackage{german,epic}
\newtheorem{Satz}{Satz}[section]
\newtheorem{Definition}[Satz]{Definition}
\newcommand{\HH}{{\bf H}} 
\newcommand{\Z}{{\bf Z}} 
\newcommand{\Q}{{\bf Q}} 
\newcommand{\F}{{\bf F}} 
\newcommand{\R}{{\bf R}} 
\newcommand{\C}{{\bf C}} 
\newcommand{\ggT}{\mbox{\rm ggT}} 
\newcommand{\Harm}{\mbox{\rm Harm}} 
\newcommand{\Min}{\mbox{\rm Min}} 
\newcommand{\Aut}{\mbox{\rm Aut}} 
\newcommand{\SL}{\mbox{\rm SL}} 
\newcommand{\vol}{\mbox{\rm vol}} 
\newcommand{\strich}[1]{| \! _{\mbox{{\raisebox{-1.5ex}{{$\scriptstyle{#1}$}}}}}  }

\title{Gitter und Modulformen}

\author{Gabriele Nebe, Ulm}
\date{}
\begin{document}

\maketitle

\noindent
{\bf Abstract} \\

\\

\noindent
A main goal in lattice theory is the construction of dense lattices.
Most of the remarkable dense lattices in small dimensions have an
additional symmetry, they are modular, i.e. similar to their 
dual lattice. 
Extremal lattices are densest modular lattices,
whoses density is as high as 
the theory of modular forms allows it to be.
The theory of theta series with harmonic coefficients allows to 
classify and to construct extremal lattices as well as 
to prove that some of them are strongly perfect and hence
local maxima of the density function.

\tableofcontents

\vspace{1cm}
\noindent
Eingegangen 12.9.2002\\ 

\noindent
Gabriele Nebe, Abteilung Reine Mathematik, Universit\"at Ulm, 89069 Ulm, Germany\\
E-Mail: nebe@mathematik.uni-ulm.de, URL: www.mathematik.uni-ulm.de/ReineM/nebe

\newpage 

\section{Einleitung}

Die Theorie der Gitter hat Beziehungen zu vielen verschiedenen
Gebieten der Mathematik. 
Die meisten bekannten sch\"onen Gitter haben  
(in der Regel sogar mehrere) Konstruktionen mit Hilfe von
Codes, algebraischen Zahlringen, algebraischer Geometrie oder endlichen
Matrixgruppen, um nur einige wichtige Schlagworte zu nennen.
Umgekehrt hilft das Konzept eines Gitters in anderen Bereichen der 
Mathematik, wie z.B. in der algebraischen Topologie oder auch
in der algorithmischen 
Zahlentheorie (\cite{PohstZas}). 
In diesem Artikel m\"ochte ich \"uber einige Anwendungen der 
Theorie der Modulformen auf die Definition und Konstruktion
dichter Gitter berichten, die insbesondere in den letzten 10 Jahren
einen neuen Aufschwung in die Gittertheorie gebracht haben.
Ausgangspunkt ist die Beobachtung, dass viele der ber\"uhmtesten
 Gitter - wie z.B. das Leech Gitter in Dimension 24, das Barnes-Wall Gitter
in Dimension 16, das Coxeter-Todd Gitter in Dimension 12 - 
 eine zus\"atzliche Symmetrie haben, sie sind
{\sf modular}, also \"ahnlich zu ihrem dualen Gitter.
In \cite{queb1} verallgemeinert Quebbemann den klassischen,
mit Hilfe von Modulformen definierten
Begriff des extremalen unimodularen Gitters auf gewisse modulare Gitter
und legt so den Grundstein f\"ur die Entdeckung
neuer extremaler modularer Gitter.
Dabei liefern Modulformen nicht nur die Definition von Extremalit\"at
sondern sind auch ein Hilfsmittel, solche Gitter zu konstruieren
oder deren Nichtexistenz zu beweisen.

Ein durch die Anwendung in der Informations\"ubertragung
motiviertes Ma\ss \ f\"ur die G\"ute eines Gitters ist seine 
{\sf Dichte}, also die Dichte der zugeh\"origen gitterf\"ormigen 
Kugelpackung. 
Die sich daraus ergebende Frage, die dichtesten
Gitter in jeder Dimension zu konstruieren, ist schon ab Dimension
9 ungel\"ost. Daher beschr\"ankt man sich beim Suchen dichter Gitter auf
gewisse sch\"one Teilklassen von Gittern, wie z.B. die in Abschnitt 
\ref{MODULAR} eingef\"uhrten modularen Gitter. 
Der Bezug zu Modulformen geschieht \"uber die 
Theta-Reihe des Gitters, das ist die erzeugende 
Funktion der Anzahlen von Gittervektoren gegebener L\"ange.
Als holomorphe Funktion auf der oberen Halbebene betrachtet, hat 
die Theta-Reihe eines modularen Gitters
zus\"atzliche Invarianzeigenschaften  und ist somit 
eine Modulform zu einer recht gro\ss en Modulgruppe. 
Durch die Kenntnis des endlich dimensionalen Vektorraums dieser
Modulformen kann man die Vielfachheit des Wertes 1
bei $i\infty $ nach oben beschr\"anken.
Eine Modulform hei\ss t extremal, wenn diese
Vielfachheit maximal wird.
Ist diese Modulform die Theta-Reihe eines modularen Gitters, so
bedeutet dies, dass das Minimum und damit auch die Dichte dieses
Gitters maximal ist unter den entsprechenden modularen Gittern.
Dann hei\ss t dieses Gitter {\sf extremal}.
Extremale Gitter sind also dichteste modulare Gitter,
f\"ur die die Dichte so gro\ss  \ ist, wie es die Theorie der Modulformen
erlaubt (siehe Abschnitt \ref{EXT} f\"ur eine genaue Definition).

Extremale Gitter liefern oft gute sph\"arische Designs.
Umgekehrt kann man mit sph\"arischen Designs lokal dichteste 
Gitter finden (Abschnitt \ref{STP}) und mit Hilfe der
Theorie der Theta-Reihen mit harmonischen
Koeffizienten zum Beispiel zeigen, dass gewisse extremale Gitter 
lokale Maxima der Dichtefunktion sind (Abschnitt \ref{HARM}). 
Diese Theorie hilft auch zum Beweis der Nichtexistenz 
und, weitaus interessanter, bei der Konstruktion und Klassifikation
extremaler Gitter \cite{BaV}.

Dieser Artikel geht bewu\ss t nicht zu stark auf die technischen
Details ein. Der interessierte Leser sei dazu auf die Originalarbeiten
\cite{queb1}, \cite{queb2} sowie den sehr sch\"onen \"Ubersichtsartikel
\cite{SchaSchuPi} verwiesen.
\"Uber  den Zusammenhang mit sph\"arischen Designs informiert \cite{Mart}
(insbesondere die Ausarbeitung \cite{Venkov}).
Neben der Bibel  der Gittertheorie \cite{SPLAG} ist 
\cite{Ebeling} eine weitere sch\"one Einf\"uhrung in die behandelte 
Thematik. Alle vorkommenden Gitter findet man in der Datenbank \cite{database}.

\section{Dichte Gitter}{\label{DICHTE}}

In dieser Arbeit ist ein Gitter $L$ immer die Menge aller ganzzahligen
Linear\-kombina\-tionen von Basisvektoren $(b_1,\ldots , b_n)$
des Euklidischen Raums $(\R ^n, (,))$.
Ein durch die Anwendung von Gittern in der 
Informations\"ubertragung motiviertes
 Ma\ss \  f\"ur die G\"ute eines Gitters $L$ ist die
 {\sf Dichte} $\delta (L) $ der zugeh\"origen Kugelpackung.
Dabei werden gleichgro\ss e $n$-dimensionale Kugeln so gepackt, dass 
die Mittelpunkte der Kugeln gerade die Gitterpunkte sind.
Der maximal m\"ogliche gemeinsame Durchmesser der Kugeln ist gleich dem minimalen
Abstand verschiedener Gitterpunkte also die
Quadratwurzel aus dem {\sf Minimum} von $L$, 
$$\min (L) := \min \{ (x,x) \mid 0\neq x \in L \} .$$
Der Raum, der pro Kugel der Packung ben\"otigt wird, ist
die Wurzel aus der {\sf Determinante}  von $L$
$$\det (L) := \vol (\R^n/L) ^2 := \det ((b_i,b_j)_{i,j=1}^n) .$$
Also ist die Dichte der Kugelpackung
$$\delta (L) = \frac{V_n}{2^n} \sqrt{\frac{\min(L)^{n}}{\det (L)}} ,$$ wo 
$V_n$ das Volumen der $n$-dimensionalen Einheitssph\"are bezeichnet.

Dichte Gitter liefern gute fehlerkorrigierende Codes f\"ur analoge Signale:
Ein Signal kann man sich als Punkt im $\R^n$ vorstellen.
\"Uber einen (nicht st\"orungsfreien) Kanal kann man solche 
Signale  nicht genau \"ubertragen. Daher l\"a\ss t man (\"ahnlich wie in der Codierungstheorie)
nur Gitterpunkte als Informationssignale zu. Zum Decodieren des gest\"orten empfangenen
Signals, mu\ss \ dann der n\"achstgelegene Gitterpunkt gefunden werden.
Um dabei m\"oglichst gro\ss e Fehler korrigieren zu k\"onnen, sollen 
die Gitterpunkte weit auseinanderliegen
(d.h. das Gitter soll ein gro\ss es Minimum haben).
 Andererseits will man mit beschr\"ankter
Energie m\"oglichst viel Information \"ubertragen (d.h. das Gitter soll 
eine m\"oglichst kleine Determinante haben).
Dabei ist die Wahl des richtigen Gitters wesentlich:
So kann man z.B. mit dem {\sf Leech Gitter} $\Lambda _{24} \subset \R^{24}$,
dem wahrscheinlich dichtesten 24-dimensionalen Gitter,
bei gleicher Energie und Fehlerkorrektur $2^{24} $ (also mehr als 16  Millionen) mal so viel 
Information \"ubertragen wie mit dem Standardgitter $\Z^{24}$.
Im 80-dimensionalen kennt man zwei Gitter ($L_{80}$ und $M_{80}$, siehe Tabelle auf Seite \pageref{ext}),
f\"ur die die Verbesserung gegen\"uber dem Standardgitter $\Z ^{80}$ einen
Faktor von mehr als $10^{36}$ ausmacht.
Trotzdem ist es i.a. nicht sinnvoll so hochdimensionale Gitter einzusetzen, da man 
dann die Information nicht mehr (schnell) decodieren kann.

Ein Hauptziel der Gittertheorie ist das Finden dichter Gitter.
Dazu gibt es eine \"uber 100 Jahre alte Theorie von Korkine und Zolotareff,
die sp\"ater von Voronoi (\cite{Voronoi}) verfeinert wurde:
Die Dichte Funktion $\delta $ hat auf dem Raum der \"Ahnlichkeitsklassen 
von $n$-dimensionalen Gittern nur endlich viele lokale Maxima,
sogenannte {\sf extreme Gitter}.
Die extremen Gitter kann man mit Hilfe der Geometrie ihrer k\"urzesten
Vektoren 
$$\Min (L) :=\{ x\in L \mid (x,x) = \min (L) \} $$
charakterisieren: 

\begin{Satz} (\cite{Voronoi})  {\label{eutper}}
$L$ ist extrem $\Leftrightarrow $ $L$ ist {\sf perfekt}
(d.h. $ \langle p_x := x^{tr} x \mid x\in \Min(L) \rangle _{\R } = \R ^{n\times n}_{sym} $) und $L$ ist {\sf eutaktisch}
(d.h. $I_n = \sum _{x\in \Min (L)} \lambda _x p_x $ mit $\lambda _x >0 \ \forall x $).
\end{Satz}

Wohingegen Eutaxie eine Konvexit\"atsbedingung ist (die die Bedingung $\delta '' (L) < 0$ ersetzt, da $\delta $ im allgemeinen nicht differenzierbar ist)
ist Perfektion eine lineare Bedingung (stellvertretend f\"ur $\delta '(L) = 0$).

Die Voronoische Charakterisierung erlaubt es u.a. zu zeigen, dass 
extreme Gitter \"ahnlich sind zu {\sf ganzen Gittern}, d.h. die 
Bilinearform $(,)$ nimmt auf $L$ nur ganze Werte an, oder gleichbedeutend
$L$ liegt in seinem {\sf dualen Gitter} 
$$L^* := \{ v\in \R^n \mid (v,x)\in \Z \ \forall x\in L \} .$$
Es ist also keine Einschr\"ankung, sich bei der Suche nach dichtesten 
Gittern auf ganze Gitter zu beschr\"anken.

Voronoi gab 1908 einen Algorithmus an, 
um alle (endlich vielen) \"Ahnlichkeitsklassen
 perfekter Gitter der Dimension $n$ aufzulisten, der bis Dimension 7 
praktikabel ist. 
Die dichtesten Gitter sind bis zur Dimension 8 bekannt, wobei es in
Dimension 8 nicht mehr m\"oglich ist, alle extremen Gitter zu bestimmen
(man kennt schon mehr als 10 000). 
Die dichtesten Gitter der Dimension $\leq 8$ sind sogenannte 
{\sf Wurzelgitter}, das sind ganze Gitter, die von Vektoren der L\"ange 2
erzeugt werden. 

\begin{center}
{\bf Die dichtesten Gitter in Dimension $n \leq 8$:}
\end{center}
\begin{center}
\begin{tabular}{|c|c|c|c|c|c|c|c|c|}
\hline
n & 1 & 2 & 3 & 4 & 5 & 6 & 7 & 8 \\
\hline
& $\Z $ & $A_2$ & $A_3$ & $D_4 $ & $D_5$ & $E_6$ & $E_7$ & $E_8 $ \\
\hline
\end{tabular}
\end{center}

Die Wurzelgitter sind alle klassifiziert und gut untersucht (siehe z.B.
\cite[Section 1.4]{Ebeling}, \cite[Chapter 4]{SPLAG}). 
Sie spielen in vielen Gebieten der Mathematik eine wichtige Rolle 
wie z.B. in der Theorie der Liealgebren und 
der algebraischen Gruppen.
In Dimension 9 ist das dichteste bekannte Gitter, das geschichtete
Gitter $\Lambda _9$ (\cite{SPLAG}), dichter als alle Wurzelgitter.

\section{Modulare Gitter}{\label{MODULAR}}

Der Begriff des modularen Gitters, der in dem hier verwendeten Sinn vor ca. 10
Jahren von H.-G. Quebbemann gepr\"agt wurde, kommt zum einen daher, dass die Theorie der 
Modulformen hilft, modulare Gitter zu untersuchen, zum anderen, da diese 
Gitter eine Verallgemeinerung der {\sf unimodularen} Gitter sind, das sind
ganze Gitter der Determinante 1.

\begin{Definition}
Ein ganzes Gitter $L$ der Dimension $2k$ 
hei\ss t {\sf $N$-modular} oder
{\sf modular der Stufe $N$}, falls $L$ isometrisch zu seinem reskalierten dualen Gitter
$$L _{(N)} :=\sqrt{N} L^* $$
ist, d.h.
 es gibt eine orthogonale Abbildung $\sigma \in O_n(\R)$ mit
$ \sigma L  = L _{(N)} $.
1-modulare Gitter hei\ss en auch {\sf unimodular}.
\end{Definition}

Viele der dichtesten Gitter in kleinen Dimensionen sind modular:
Die Wurzelgitter $A_2$, $D_4$ und $E_8$, welche die
dichtesten Gitter in Dimension 2, 4 und 8 sind, sind $3$-, $2$- bzw. unimodular.
Weitere ber\"uhmte modulare Gitter sind das 
$3$-modulare Coxeter-Todd Gitter $K_{12}$ in Dimension 12,
das $2$-modulare Barnes-Wall Gitter $BW_{16}$ in Dimension 16 
und nat\"urlich das bemerkenswerte Leech Gitter $\Lambda _{24}$,
das einzige gerade unimodulare Gitter der Dimension 24, welches keine
Vektoren der L\"ange 2 besitzt.
In all diesen F\"allen ist die Stufe $N$ eine Primzahl (oder 1). 
W\"ahrend der Klassifikation der maximal endlichen rationalen Matrixgruppen
(\cite{NeP}, \cite{dim24}, \cite{dim2531}) wurden einige sehr dichte 
modulare Gitter gefunden, f\"ur welche die Stufe
 eine quadratfreie zusammengesetzte Zahl ist.
Ist  $N=mm'$ ein Produkt von 2 teilerfremden Zahlen, 
so gibt es zwischen dem
$N$-modularen Gitter $L$ und seinem dualen Gitter $L^*$ zwei weitere sogenannte
``partielle'' duale Gitter $$L^{*,m}:= \frac{1}{m} L \cap L^*$$ und $L^{*,m'}$. 

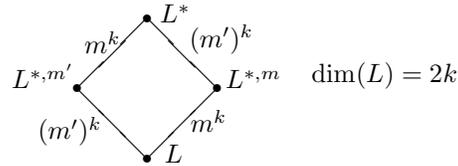
$$
\setlength{\unitlength}{0.0003in}
\begingroup\makeatletter\ifx\SetFigFont\undefined%
\gdef\SetFigFont#1#2#3#4#5{%
  \reset@font\fontsize{#1}{#2pt}%
  \fontfamily{#3}\fontseries{#4}\fontshape{#5}%
  \selectfont}%
\fi\endgroup%
{\renewcommand{\dashlinestretch}{30}
\begin{picture}(5956,2400)(0,-10)
\put(2074,83){\circle*{150}}
\put(3274,1283){\circle*{150}}
\put(2074,2483){\circle*{150}}
\put(874,1283){\circle*{150}}
\drawline(2074,2483)(874,1283)
\drawline(2074,2483)(3274,1283)
\drawline(874,1283)(2074,83)
\drawline(3274,1283)(2074,83)
\put(2374,8){\makebox(0,0)[lb]{$L$}}
\put(2299,2408){\makebox(0,0)[lb]{$L^*$}}
\put(1000,1883){\makebox(0,0)[lb]{$m^k$}}
\put(2824,533){\makebox(0,0)[lb]{$m^k$}}
\put(2800,1883){\makebox(0,0)[lb]{$(m')^k$}}
\put(1300,300){\makebox(0,0)[rb]{$(m')^k$}}
\put(4900,1283){\makebox(0,0)[lb]{$\dim (L) =2k$}}
\put(824,1283){\makebox(0,0)[rb]{$L^{*,m'}$}}
\put(3424,1283){\makebox(0,0)[lb]{$L^{*,m}$}}
\end{picture}
}
$$

\begin{Definition}
Ein $N$-modulares Gitter $L$ hei\ss t {\sf stark $N$-modular}, falls $L$ 
isometrisch ist zu allen reskalierten partiellen dualen Gittern 
$$L_{(m)} := \sqrt{m} L^{*,m} $$ f\"ur alle exakten 
Teiler $m$ von $N$ (d.h. $\ggT (m, N/m) = 1$).
\end{Definition}

Da es zu schwierig ist, die absolut dichtesten Gitter zu bestimmen, und (zumindest in
``kleinen'' Dimensionen) viele der dichtesten bekannten Gitter modular sind, ist es 
ein interessantes Problem, die dichtesten (stark) $N$-modularen Gitter zu finden.
Da die Determinante eines $N$-modularen Gitters $L$ der Dimension $2k$
gleich $\det (L) = N^k $ ist, bedeutet dies, die $N$-modularen Gitter
mit dem gr\"o\ss tm\"oglichen Minimum zu bestimmen. 
Dabei hilft die Theorie der Modulformen, dieses Minimum nach oben zu beschr\"anken, 
so dass man einem einzelnen stark $N$-modularen Gitter ansieht, ob es ein dichtestes
modulares Gitter ist, ohne alle anderen Gitter zu kennen.

\section{Modulformen}

Die im n\"achsten Abschnitt definierten 
Theta-Reihen von Gittern haben gewisse Invarianzeigenschaften unter
Variablensubstitutionen, sie sind Modulformen. 
Da diese Beobachtung grundlegend f\"ur diesen Artikel ist, wird 
kurz auf den Begriff der Modulform eingegangen.
Wie man in der Funktionentheorie lernt, 
ist die Gruppe der  biholomorphen Abbildungen 
auf der oberen Halbebene $\HH := \{ z\in \C \mid \Im (z) > 0 \} $
die Gruppe der M\"obiustransformationen 
$$z \mapsto A(z) := \frac{az+b}{cz+d} , \ \ A = \left( \begin{array}{cc} a & b \\ c & d \end{array} \right) \in SL_2(\R ) .$$
Dies liefert f\"ur alle $k\in \Z$ eine Operation $\strich{k}$
der $\SL_2(\R )$ auf dem Vektorraum der meromorphen
Funktionen $f: \HH \rightarrow \C $ 
durch 
$$ f \strich{k}   A  (z) := (cz+d)^{-k} f (\frac{az+b}{cz+d} ) .$$

\begin{Definition}
Sei $U $ eine  Untergruppe von $SL_2(\R )$, so dass 
$U\cap SL_2(\Z) $ in $U$ und $SL_2(\Z )$ endlichen Index hat, 
und $\chi : U \rightarrow \C^*$ ein
Charakter. Eine meromorphe Funktion $f: \HH \rightarrow \C $ hei\ss t
eine {\sf Modulform vom Gewicht $k$ zum Charakter $\chi $},
$f\in {\cal M}_k(U,\chi )$, falls 
$$ f \strich{k} A  = \chi (A) f \mbox{ f\"ur alle } A\in U $$
gilt und $f\strich{k} M$ f\"ur alle $M \in SL_2(\Z)$ bei $i \infty $ h\"ochstens einen Pol hat.
$f$ hei\ss t eine {\sf Spitzenform},
$f\in {\cal S}_k(U,\chi )$,  
falls zus\"atzlich $\lim _{t\to \infty} f \strich{k} M (it ) = 0 $ f\"ur alle $M\in SL_2(\Z )$.
\end{Definition}

Da die $\strich{k} $-Operation  ``multiplikativ'' ist, d.h.
$(f \strich{k}  A )(g \strich{l}  A )= (fg) \strich{k+l} A $ ist 
f\"ur jede Familie von Charakteren
$\chi _k: U \rightarrow \C ^*$, $k\in \Z _{\geq 0}$ mit 
$\chi _k \chi _l = \chi _{k+l}$ der Raum 
$${\cal M}(U, (\chi _k)) = \bigoplus _{k=0}^{\infty } {\cal M}_k (U,\chi _k)$$  ein graduierter Ring,
der {\sf Ring der  Modulformen}
von $U$
 (zum Charakter $(\chi _k)$), in
welchem die Spitzenformen ein Ideal bilden.
Dieser graduierte Ring ist endlich erzeugt.
Die Dimensionen der  graduierten Komponenten lassen sich mit Hilfe von
Spurformeln berechnen (siehe z.B. \cite{Miyake}).
Die Modulformen $f$, die f\"ur ganze Gitter
von Bedeutung sind, sind alle invariant unter der Transformation $z\mapsto z+2$.
Also haben sie eine Fourier-Entwicklung 
$$f(z) = \sum _{n=0}^{\infty } a_n q^n , \ \ q = \exp(\pi i z ) .$$
Kennt man die Dimension von 
${\cal M}_k(U,\chi _k)$, so gen\"ugt es 
f\"ur festes $k$ mit einer festen Genauigkeit (also modulo $q^a$ f\"ur hinreichend gro\ss es $a$) zu rechnen.
Aus der Potenzreihenentwicklung der 
Ringerzeuger modulo $q^a$ erh\"alt man durch Multiplikation dann eine explizite 
Basis von
${\cal M}_k(U,\chi _k)$.

\section{Theta-Reihen}{\label{THETA}}
Sei $L$ ein ganzes Gitter. Die {\sf Theta-Reihe} $\theta _L$ von 
$L$ ist die erzeugende Funktion der Anzahlen von Vektoren gegebener L\"ange 
in $L$, 
$$\theta _L (z) := \sum _{j=0}^{\infty } a_L(j) q^j $$
wo $a_L(j) = |\{ x\in L \mid (x,x) = j \} |$ die Anzahl der Vektoren der L\"ange $j$ in
$L$ ist und $$q:=\exp (\pi i z).$$
Dann ist $\theta _L$ eine holomorphe Funktion auf der oberen Halbebene.
Da $\exp (2\pi i ) = 1$ ist, ist z.B. $\theta _L(z) = \theta _L(z+2)$ und sogar 
$\theta _L (z) = \theta _L(z+1)$, falls $L$ ein {\sf gerades Gitter} ist,
d.h. $(x,x) \in 2\Z $ f\"ur alle $x\in L$.
Sei nun $L$ ein gerades Gitter.
Dann ist die {\sf Stufe} $N$ von $L$ die kleinste nat\"urliche Zahl, f\"ur die
$$L_{(N)} := \sqrt{N} L^*$$ ein gerades Gitter ist.

\begin{Satz}
Ist $L$ ein gerades Gitter der Dimension $2k$ von Stufe $N$, so ist 
$$\theta _L \in {\cal M}_k (\Gamma _0 (N) , \chi _k ) $$
eine Modulform zur Gruppe 
$$\Gamma _0 (N) : = \{ \left( \begin{array}{cc} a & b \\ c & d \end{array} \right) 
\in \SL _2(\Z ) \mid  N \mbox{ teilt } c \} $$ 
mit dem Charakter $\chi _k$ definiert durch das Legendre Symbol
$$\chi _k 
(\left( \begin{array}{cc} a & b \\ c & d \end{array} \right) ) := 
(\frac{(-1)^k \det(L)}{d}) \in \{ \pm 1 \} .$$
\end{Satz}

Im Hinblick auf modulare Gitter
interessieren gerade Gitter $L$ der Dimension $2k$ von
Stufe $N$ der Determinante $N^k$.
Diese Gitter bilden f\"ur Primzahlen $N$ oder $N=1$ ein {\sf Geschlecht},
d.h. f\"ur je zwei solche Gitter
$L$ und $M$ sind ihre  Lokalisierungen 
$\Z_p\otimes M \cong \Z _p \otimes L$ 
f\"ur alle Primzahlen $p$ isometrisch.

Die meisten Standardkonstruktionen f\"ur Gitter haben eine Entsprechung auf der Ebene 
der Theta-Reihen: 
Mit Hilfe von Poisson Summation erh\"alt man  die 
Theta-Transformationsformel 
$$ 
\theta _{L^*} (z) =
 (\frac{z}{i} )^{-k} \sqrt{\det (L) } \theta _{L}(-\frac{1}{z} ) 
 \ \ \ \  (\mbox{wo } 2k=\dim(L))$$
welche die Theta-Reihe des dualen Gitters durch  $\theta _L$ ausdr\"uckt.

Ist $\det (L) = N^{k}$ so ergibt sich
nach Substitution von $z$ durch $\sqrt{N} z $ in der Theta-Transformationsformel 
dass 
$$\theta _{L_{(N)}}  = \chi _{N,k} (t_N) \theta _{L} \strich{k} t_N  \ \ \ \ (\mbox{wo } L_{(N)} = \sqrt{N} L^*) $$ 
f\"ur die
Fricke Involution $$ 
 t_N :=\left( \begin{array}{cc} 0 & \sqrt{N}^{-1} \\ -\sqrt{N} & 0 \end{array} \right)
\in N_{SL_2(\R)} (\Gamma _0(N)) $$ gilt, wobei 
$\chi _{N,k} (t_N ) := i^k $ gesetzt ist.

Etwas allgemeiner, kann man auch die Theta-Reihe der partiellen dualen 
Gitter aus $\theta _L$ berechnen (\cite[p. 60]{queb2})
und findet, dass die Theta-Reihe eines stark $N$-modularen Gitters eine
Eigenfunktion der Atkin-Lehner Involution $W_m \in N_{SL_2(\R)}(\Gamma _0(N))$
 f\"ur alle exakten Teiler $m$ von $N$ ist.

\begin{Satz}
Ist $L$ ein gerades stark $N$-modulares Gitter der Dimension $2k$, so ist $\theta _L$ eine
Modulform zur Gruppe
$$\Gamma _*(N) := \langle 
\Gamma _0(N),  W_m  \mid   m \mbox{ exakter Teiler von } N 
\rangle  $$
zu einem Charakter $\chi _{N,k}$ mit $\chi _{N,k}\chi _{N,k'} = \chi _{N,k+k'}$.
\end{Satz}

Diese Atkin-Lehner Involutionen erzeugen eine elementar abelsche $2$-Gruppe
$W(N) \leq N_{SL_2(\R)}(\Gamma _0(N)) / \Gamma _0 (N) $ deren
Rang  die Anzahl der Primteiler von $N$ ist.
$W(N)$ operiert auf dem $\C$-Vektorraum, dessen Basis die Isometrieklassen 
der geraden Gitter von Stufe $N$, Dimension $2k$ und Determinante $N^k$ bilden, durch
$$[L] \cdot W_m := \chi _{N,k}(W_m)  [L_{(m)} ].$$
Unter dem linearen Operator ``Theta-Reihe nehmen'' 
wird diese Operation gerade auf die Operation von $W(N)$ auf dem Raum der
Modulformen f\"ur $\Gamma _0(N) $ abgebildet.
Der Vorteil dieser schon bei Eichler \cite[Kapitel IV]{Eichler}
beschriebenen Sichtweise
 liegt vor allem darin, dass man auch andere
Operatoren z.B. ``Theta-Reihe mit sph\"arischen Koeffizienten nehmen''
oder ``Siegelsche Theta-Reihe bilden'' anwenden kann und dann die
entsprechende Operation auf den (Siegelschen) Modulformen erh\"alt.
Diese Philosophie wurde z.B. in \cite{Siegel} verfolgt.

\section{Extremale Gitter}{\label{EXT}}
Der Ring der Modulformen von $\Gamma _*(N)$ ist besonders \"ubersichtlich, 
wenn die Summe $\sigma _1(N)$ der Teiler von $N$ ein Teiler von 24 ist.
Daher wird im ganzen Abschnitt vorausgesetzt, dass 
 $$ N \in \{ 1,2,3,5,6,7,11,14,15,23 \} =: {\cal A} $$
ist.
Sei weiter $L$ ein gerades stark $N$-modulares Gitters minimaler
Dimension $2d_N$ und
$\theta _N $ seine Theta-Reihe (vom Gewicht $d_N$).
Dann liegen die Theta-Reihen der
stark $N$-modularen Gitter, welche zum Geschlecht von $L^d$ f\"ur ein $d$ 
geh\"oren, in dem  Ring 
$${\cal M}(N) := 
  \C [\theta _N , \Delta _N ] \subset
\bigoplus _{k=0}^{\infty } {\cal M}_{k}(\Gamma _*(N) ,\chi _{N,k} ) $$
wo
$$\Delta _N := \prod _{m \mid N} \eta (mz) ^{24/\sigma_1(N)}$$ 
eine Spitzenform vom Gewicht $k_N := \frac{12\sigma_0(N)}{\sigma _1(N)}$ ist.
Hier ist $$\eta (z) = q^{1/12} \prod _{n=1}^{\infty } (1-q^{2n}) $$
 die Dedekindsche $\eta $-Funktion
und $\sigma _0(N)$ die Anzahl der Teiler von $N$.
Die Reihe $\theta _N$ und der Charakter $\chi _{N,k}$ und damit auch
${\cal M}(N)$ h\"angen
vom Geschlecht des gew\"ahlten Gitters $L$  ab (welches jedoch f\"ur $N\neq 6$
durch $d_N$ eindeutig bestimmt ist),
$\Delta _N$ ist schon durch die Stufe  $N$ bestimmt.
$\Delta _N$ und $\theta _N$ sind Potenzreihen in $q^2$.
Die $q$-Entwicklung von $\Delta _N$ beginnt mit $q^2$ und die von $\theta _N$ mit 1.
$d_N$ und $k_N$ ergeben sich aus der folgenden Tabelle:
$$\begin{array}{|c|c|c|c|c|c|c|c|c|c|c|}
\hline
N & 1 & 2 & 3 & 5 & 6 & 7 & 11 & 14 & 15 & 23 \\
\hline
2k_N & 24 & 16 & 12 & 8 & 8 & 6 & 4 & 4 & 4  & 2\\
\hline
2d_N & 8 & 4 & 2 & 4 & 4 & 2 & 2 & 4 & 4  & 2 \\
\hline
\end{array}
$$

Ist $L$ ein gerades stark $N$-modulares Gitter, so gibt es also $a_i \in \C$
($a_0 = 1$) mit
$$\theta _{L} = 
\sum _{i=0}^l a_i \Delta _N ^i \theta _N ^j ,\mbox{ wo }  k_N i + d_N j = k 
\mbox{ und } l = \lfloor \frac{k}{k_N} \rfloor $$ die gr\"o\ss te ganze Zahl 
$\leq \frac{k}{k_N} $ bezeichnet.
Die $q$-Entwicklung von $\Delta _N^i \theta _N^j$ beginnt mit $q^{2i}$. 
Also enth\"alt ${\cal M}(N)$  genau eine Funktion  
$f_{N,k} $ vom Grad $k$ mit $q$-Entwicklung 
$$f_{N,k} = 1+0\cdot q^2 + \ldots + 0 \cdot q^{2l} + a_{2l+2} \cdot q^{2l+2} + \ldots \in {\cal M} _{k}(N)$$
die sogenannte {\sf extremale Modulform} vom Gewicht $k$.
Schon C.L. Siegel \cite{CLSiegel} hat f\"ur $N=1$ gezeigt, 
dass der Koeffizient $a_{2l+2} $ echt positiv ist
(siehe \cite[Theorem 2.0.1]{SchaSchuPi} f\"ur $N>1$).
Also ist das Minimum eines geraden 
stark $N$-modularen Gitters  der Dimension $2k$ immer kleiner oder gleich
$2+2\lfloor \frac{k}{k_N}\rfloor $.

\begin{Definition}
Sei $N\in {\cal A} $ und $k_N := \frac{12\sigma_0(N)}{\sigma _1(N)}$
das Gewicht der Spitzenform $\Delta _N$.
Ein gerades stark $N$-modulares Gitter $L$ der Dimension $2k$ hei\ss t 
{\sf extremal}, falls $$\min (L) = 2+2\lfloor \frac{k}{k_N} \rfloor  .$$
\end{Definition}

Grob gesprochen bedeutet dies, dass ein Gitter extremal ist, wenn sein Minimum
so gro\ss \ ist, wie es die Theorie der Modulformen erlaubt.
Mit dieser Philosophie kann man auch Extremalit\"at f\"ur andere Geschlechter
von Gittern definieren (siehe \cite[Definition 1.6]{SchaSchuPi}).
Ist $L$ also ein extremales stark $N$-modulares Gitter der Dimension $2k$,
 so ist $\theta _L = f_{N,k}$ 
die (dem Geschlecht von $L$ zugeordnete) extremale Modulform.
Eine notwendige Bedingung f\"ur die Existenz eines extremalen geraden stark $N$-modularen
Gitters in Dimension $2k$ ist, dass die Koeffizienten $a_{m}$ der
extremalen Modulform $f_{N, k}$ f\"ur $m>0$ nichtnegative gerade Zahlen sind.
In \cite{MOS} wird gezeigt, dass f\"ur gro\ss e $k$ $(k>20500)$
der Koeffizient $a_{2l+4}$ von $f_{1,k}$ negativ wird.
In Dimensionen $\geq 41000 $ gibt es also keine extremalen geraden unimodularen
Gitter mehr.  
Jedoch ist schon in Dimension 72
die Existenz eines extremalen Gitters ein immer noch offenes Problem.
Mit derselben Technik zeigt \cite[Theorem 2.0.1]{SchaSchuPi} das analoge Resultat
f\"ur $N>1$. 
Es gibt also f\"ur festes $N\in {\cal A}$
nur endlich viele extremale gerade stark $N$-modulare Gitter.

Da die Definition der  Extremalit\"at nur die Theta-Reihen der Gitter benutzt,
gelten die Schranken an die Minima auch f\"ur ``formal'' stark $N$-modulare
Gitter, also solche Gitter, deren Theta-Reihe gleich der Theta-Reihe
aller partiellen dualen Gitter ist.

Die Dimension
$2k_N$, also die erste Dimension, in der ein extremales gerades
stark $N$-modulares Gitter Minimum 4 hat, ist besonders bemerkenswert.
Es gibt n\"amlich je genau ein extremales stark $N$-modulares
Gitter $E^{(N)}$
der Dimension $2k_N$.
F\"ur $N=1$ ist dieses Gitter das Leech Gitter $\Lambda _{24}$.
 Eine einheitliche Konstruktion dieser extremalen Gitter ist in
\cite{RaS} gegeben:
 Die Zahlen $N$ sind genau die Elementordnungen in der Matthieu Gruppe
$M_{23} \leq \Aut(\Lambda _{24})$, die auf dem Leech Gitter  als
Automorphismen operiert. 
Dann ist $E^{(N)}$ das Fixgitter in $\Lambda _{24}$ 
eines Elements der Ordnung $N$ in $M_{23}$.

Die Vielfachen von $2k_N$ hei\ss en {\sf Sprungdimensionen} f\"ur die
stark $N$-modularen Gitter. Extremale Gitter in diesen Dimensionen sind von
besonderem Interesse, da sie meist sehr dicht sind. 
In der Regel gibt es in diesen Dimensionen nur sehr wenige extremale Gitter,
so kennt man z.B. nur 4 extremale gerade unimodulare Gitter  in Sprungdimensionen
(siehe folgende Tabelle).

\newpage

\begin{center} {\bf Ausgew\"ahlte extremale Gitter \label{ext}
\footnote{Eine ausf\"uhrliche und 
st\"andig aktualisierte Version dieser Tabelle und aller vorkommenden
Gitter findet man unter \cite{database}.}
.}
\end{center}
\begin{center}
\begin{tabular}{|c|c|c|c|l|}
\hline
Stufe 
& dim & min & Anzahl & Gitter  \\
\hline
N=1 & 8 & 2 & 1 & $E_8$  \\
\cline{2-5}
& 16 & 2 & 2 & $E_8\perp E_8$, $D_{16}^+$  \\
\cline{2-5}
& 24 & 4 & 1 & $\Lambda _{24}$  (Leech) \\
\cline{2-5}
& 32 & 4 & \multicolumn{2}{l|}{$\geq 10$ Millionen \cite{King}} \\
\cline{2-5}
& 48 & 6 & $\geq 3$ & $P_{48p}$, $P_{48q}$ \cite[p. 195]{SPLAG}, $P_{48n}$
\cite{cyclo} \\
\cline{2-5}
& 72 & 8 & ?  & \\
\cline{2-5}
& 80 & 8 & $\geq 2$ & $L_{80}$, $M_{80}$ \cite{BaN} \\
\hline 
N=2 & 4 & 2 & 1 & $D_4$  \\
\cline{2-5} 
& 8 & 2 & 1 & $D_4\perp D_4$  \\
\cline{2-5} 
& 12 & 2 & 3 \cite{SchaHe} &  \\ 
\cline{2-5} 
& 16 & 4 & 1 \cite{SchV} & $BW_{16}$ \cite{BaW}  \\
\cline{2-5} 
& 20 & 4 & 3 \cite{BaV} & $[(SU_5(2) \circ SL_2(3)).2]_{20}$, $[2.M_{12}.2]_{20}$ \cite{NeP} \\
& & & &  $HS_{20}$  \footnotemark\ \\
\cline{2-5} 
& 32 & 6 & $\geq 4$ & $Q_{32} $ \cite{queb32}, $Q_{32}'$ \cite{queb32q} ,
$B_{32}$ \cite{Bachochurw} \\
& & & & $CQ_{32}$ \cite[Theorem 5.1]{cyclo}  \\
\cline{2-5} 
& 48 & 8 & \multicolumn{2}{l|}{$\geq 2$ \cite[Theorem 6.7]{Bachoc} } \\
\hline
N=3 & 2 & 2 & 1 & $A_2$ \\
\cline{2-5}
& 12 & 4 & 1 \cite{SchV2} & $K_{12}$ (Coxeter-Todd) \\
\cline{2-5}
& 14 & 4 & $\geq 1$  & $[\pm G_2(3)]_{14}$ \cite{KKU} \\
\cline{2-5}
& 24 & 6 & $\geq 1 $ & $[SL_2(13) \circ SL_2(3)]_{24}$ \cite{dim24} \\
\cline{2-5}
& 26 & 6 & $\geq 1 $ & $[(\pm S_6(3) \times C_3 ).2]_{26}$ \cite{dim2531} \\
\cline{2-5}
& 36 & 8 & ? & \\
\cline{2-5}
& 40 & 8 & $\geq 1$ & $L_{40}$ \cite{BaN} \\
\cline{2-5}
& 64 & 12 & $\geq 1$ &  $^{(p3)} L_{8,2} \otimes_{\infty , 3} L_{32,2}$ \\
& & & & \cite[Remark 5.2]{cyclo}, 
\cite[Prop. 4.3]{seoul}
 \\
\hline
N=5 & 8 & 4 & 1 & $H_4$ \cite{ShT} \\
\cline{2-5} 
& 16 & 6 & 1 \cite{BaV} & $[2.Alt_{10}]_{16}$  \cite{NeP}  \\
\cline{2-5}
& 24 & 8 & $\geq 1$ & $[2.J_2\circ SL_2(5) .2 ]_{24}$ \cite{Tits} \\
\hline
N=6 
 & 8 & 4 & 1 & $A_2\otimes D_4$  \\
\cline{2-5}
& 16 & 6 & $\geq 1 $ & $[((Sp_4(3) \circ C_3 )\otimes_{\sqrt{-3}} SL_2(3) ).2 ]_{16}$ \cite{NeP} \\
\cline{2-5}
& 24 & 8 & $\geq 2 $ & 
  $[(SL_2(3)\circ C_4).2 \otimes_{\sqrt{-1}} U_3(3) ]_{24}$ 
\\ & & & &
$[(6.L_3(4).2\otimes D_8 ).2 ]_{24}$
\cite{dim24} \\
\hline
N=7 & 6 & 4 & 1 & $A_6^{(2)}$  \\
\cline{2-5}
 & 12 & 6 &  0 \cite{SchaHe} &  \\
\cline{2-5}
 & 18 & 8 &  0 \cite{BaV} &  \\
\hline
N=11 & 4 & 4 & 1 &  \\
\cline{2-5}
 & 8 & 6 & 1 \cite{SchaHe} & \cite[Theorem 5.1]{cyclo} \\
\cline{2-5}
 & 12 & 8 & 0 \cite{NeV} & \\
\hline
\end{tabular}
\end{center}
\footnotetext{R. Scharlau, B. Hemkemeier durch Suche im Nachbarschaftsgraph}  
\newpage
\begin{center}
\begin{tabular}{|c|c|c|c|l|}
\hline
Stufe 
& dim & min & Anzahl & Gitter  \\
\hline
N=14 & 4 & 4 & 1 &  \\
\cline{2-5}
 & 8 & 6 & $1$ \cite{SchaSchuPi} &  \\
\cline{2-5}
 & 12 & 8 & $\geq 1$ & $[(L_2(7) \otimes D_8).2] _{12}$ \cite{NeP} \\
\hline
N=15 & 4 & 4 & 1 &  \\
\cline{2-5}
 & 8 & 6 & $2$ \cite{SchaSchuPi} &  \\
\cline{2-5}
 & 12 & 8 & $\geq 1$ & $A_2\otimes M_{6,2}$ \cite{NeP} \\
\cline{2-5}
 & 16 & 10 & $\geq 1$ & $[(SL_2(5) \otimes _{\infty , 3} SL_2(9) ).2]_{16} $ \cite{NeP} \\
\hline
\end{tabular}
\end{center}

Die extremalen Gitter in obiger Tabelle wurden auf unterschiedliche Weisen gefunden
und auch f\"ur den Beweis der Nichtexistenz extremaler Gitter gibt es
verschiedene Methoden.

R. Scharlau und B. Hemkemeier (\cite{SchaHe}) 
listen mit der Kneserschen Nachbarschaftsmethode
 alle Gitter in kleinen
Geschlechtern auf und finden so die extremalen stark modularen Gitter in
diesem Geschlecht oder zeigen   deren Nichtexistenz.
Viele extremale Gitter wurden von C. Bachoc mit Codes \"uber Zahlk\"orpern oder
Quaternionenalgebren konstruiert (\cite{Bachoc}, \cite{Bachochurw}, \cite{BaN}).
Eine weitere reiche Quelle f\"ur extremale Gitter ist die Klassifikation der
maximal endlichen rationalen (und quaternionialen)
Matrixgruppen \cite{NeP}, \cite{dim24}, \cite{dim2531},
\cite{quat} (siehe auch \cite{cyclo}).
Die sicherlich interessanteste Methode ist aber eine
gezielte Konstruktion oder ein Nichtexistenzbeweis mit Hilfe von Modulformen.
In \cite{NeV} wird die Nichtexistenz gewisser extremaler
Gitter mit Siegelschen Modulformen gezeigt.
Das wichtigste Hilfmittel 
 sind jedoch Modulformen mit sph\"arischen
Koeffizienten (siehe z.B. \cite{BaV}). 
Diese Theorie erlaubt es auch zu zeigen, dass f\"ur gewisse $N$ und $k$ 
extremale stark $N$-modulare Gitter der Dimension $2k$ lokale
Maxima der Dichte Funktion $\delta $
sind, wie in den folgenden beiden Abschnitten beschrieben wird.

\section{Sph\"arische Designs und stark perfekte Gitter}{\label{STP}}
In diesem Abschnitt wird auf neuere Entwicklungen in der Gittertheorie 
eingegangen, die haupts\"achlich von B. Venkov initiiert wurden.
Eine sch\"one Einf\"uhrung ist die Ausarbeitung \cite{Venkov}.

\begin{Definition}
Eine nicht leere endliche Teilmenge $X$ der $(n-1)$-dimen\-sio\-nalen Sph\"are
$S^{n-1} := \{ x\in \R^n \mid (x,x) = 1 \}$ hei\ss t (sph\"arisches) 
{\sf $t$-Design}, falls 
der Mittelwert \"uber $X$ gleich dem $O_n(\R )$-invarianten Integral  ist
$$ (\star) \ \ \frac{1}{|X|} \sum _{x\in X} f(x) = \int _{S^{n-1}} f(x) d\mu (x) $$
f\"ur alle homogenen Polynome $f\in \R [x_1,\ldots , x_n]_{k}$ vom Grad $k\leq t$.
\end{Definition}

Auf der rechten Seite von $(\star )$ steht das $O_n(\R)$-invariante Skalarprodukt
von $f$ mit der konstanten Funktion 1. 
Da die harmonischen Polynome 
$$\Harm _k (n) := \{ f\in \R [x_1,\ldots , x_n] _{k} \mid \Delta f = 0 \} $$
(wo $\Delta = \sum _{i=1}^n \frac{\partial ^2}{\partial x_i^2 }$ der {\sf Laplace
Operator} ist)
vom homogenen Grad $k$ f\"ur $k\geq 1$ senkrecht auf den konstanten Funktionen
stehen, folgt leicht, dass $X$ ein sph\"arisches $t$-Design ist, 
genau dann wenn 
$$\sum _{x\in X} f(x) = 0 \mbox{ f\"ur alle } f\in \Harm _k(n) , \ 1\leq k \leq t .$$
 
In der Anwendung auf Gitter, wird $X$ aus allen Gittervektoren gegebener 
L\"ange $a$  in $L$ (einer sogenannten {\sf Schicht} $L_a \subset L$) bestehen,
$$X= L_ a := \{ x\in L \mid (x,x) = a \}$$
(reskaliert, so dass $X\subset S^{n-1}$)
 meist sogar $X=\Min (L) = L_{\min(L)}$.
Diese Mengen $X$ sind {\sf symmetrisch}, d.h.  mit $x\in X$ liegt auch immer
der negative Vektor $-x $ in $X$.
Insbesondere gilt automatisch $\sum _{x\in X } f(x) = 0$ f\"ur alle homogenen
Polynome $f$ ungeraden Grades. 

\begin{Satz} (\cite[Th\'eor\`eme 3.2]{Venkov})
Sei $X \subset S^{n-1}$ endlich, symmetrisch, nicht leer. 
Dann ist $X$ ein $(2k+1)$-Design $\Leftrightarrow $
$$ \sum _{x\in X} (\alpha , x )^{2k} = |X| \frac{1\cdot 3 \cdots (2k-1)}{n(n+2)\cdots(n+2k-2)} (\alpha , \alpha )^k \mbox{ f\"ur alle } \alpha \in \R ^n .$$
\end{Satz}

Dabei gen\"ugt es vorauszusetzen, dass die linke Seite ein konstantes Vielfaches der $k$-ten 
Potenz der quadratischen Form ist. Die genaue Konstante ergibt sich dann durch 
Anwenden des Laplace Operators (nach $\alpha $). 

\begin{Definition}
Ein Gitter $L$ hei\ss t {\sf stark perfekt}, 
falls $\Min (L) $ ein sph\"arisches 4-Design ist.
\end{Definition}

\noindent
{\bf Beispiel.}
Ist $L$ ein orthogonal unzerlegbares Wurzelgitter der Dimension $n$, 
so gilt f\"ur alle $\alpha \in \R^n$
$$\sum _{x\in L_2} (x,\alpha )^2 = 2 h (\alpha , \alpha ) $$
wo $h = \frac{|L_2|}{n} $ die Coxeter Zahl ist
(siehe z.B. \cite[Proposition 1.6]{Ebeling}).
Irreduzible Wurzelsysteme bilden also 3-Designs.
\\
Die einzigen stark perfekten Wurzelgitter sind 
$A_1=\Z$, $A_2$,  $D_4$, $E_6$, $E_7$ und $E_8$ (\cite[Th\'eor\`eme 5.7]{Venkov}).
Dabei bilden nur die Wurzeln in $E_8$ ein $7$-Design. \\

F\"ur die k\"urzesten Vektoren $X:= \Min(L)$ eines stark perfekten Gitters $L$ 
gilt also
$$ (\star _4) \ \ \sum _{x\in X} (\alpha , x )^{4} = |X| m^2 \frac{1\cdot 3}{n(n+2)} (\alpha , \alpha )^2 \mbox{ f\"ur alle } \alpha \in \R ^n $$
wo $m:= \min (L)$ gesetzt wird.
Durch Anwenden des Laplace Operators $\Delta $  auf $(\star _4)$ findet man 
$$ (\star _2) \ \ \sum _{x\in X} (\alpha , x )^{2} = |X| m\frac{1}{n} (\alpha , \alpha ) \mbox{ f\"ur alle } \alpha \in \R ^n .$$
Setzt man $\alpha \in L^*$ in diese  beiden Gleichungen ein, so liefern sie
kombinatorische Bedingungen an $|X|$, $\min (L)$ und $\min (L^*)$, mit deren
Hilfe man in kleinen Dimensionen alle stark perfekten Gitter klassifizieren kann.
Man kennt alle stark perfekten Gitter der Dimension $\leq 11$:
\begin{center}
\begin{tabular}{|c|c|c|c|c|c|c|c|c|c|c|}
\hline
 1 & 2 & 3 & 4 & 5 & 6 & 7 & 8 & 9 & 10 & 11 \\
\hline
 $\Z $ & $A_2$ & $-$ & $D_4$ & $-$ & $E_6$, $E_6^*$ & $E_7$, $E_7^*$ & $E_8$ & $-$ & $K_{10}'$, $K_{10}'^{*}$  & $-$ \\
\hline
\end{tabular}
\end{center}
Die wichtigste Motivation, stark perfekte Gitter zu betrachten, ist der folgende Satz.

\begin{Satz}{\label{haupt}}
Stark perfekte Gitter sind extrem, also lokale Maxima der 
Dichte Funktion.
\end{Satz}

\underline{Beweis:}
Sei $L$ ein stark perfektes Gitter, $X:=\Min (L)$ und $m:=\min (L)$.
Nach Satz \ref{eutper} gen\"ugt es zu zeigen, dass $L$
eutaktisch und perfekt ist.
Mit $(\star _2)$ ergibt sich, dass $L$ eutaktisch ist mit
Eutaxiekoeffizienten $\lambda _x = \frac{n}{m|X|} $ f\"ur alle $x\in X$.
Es ist n\"amlich f\"ur alle $\alpha \in \R^n$ 
$$\alpha I_n \alpha ^{tr} = (\alpha , \alpha ) \stackrel{(\star _2)}{=}
\frac{n}{m |X|} \sum _{x\in X} (\alpha , x)^2 =
\frac{n}{m |X|} \sum _{x\in X} \alpha x^{tr} x  \alpha ^{tr}  $$
woraus die Gleichung f\"ur die Eutaxie folgt.
Perfektion zeigt man mit Hilfe der Gleichung $(\star _4)$ 
(siehe \cite[Th\'eor\`eme 6.4]{Venkov}).
\hfill{q.e.d.}

\begin{Satz}{\label{min}}
Sei $L$ ein stark perfektes Gitter.
Dann ist 
$$\min (L) \min (L^*)   \geq \frac{n+2}{3} .$$
\end{Satz}

\underline{Beweis:}
Sei $\alpha \in \Min (L^*)$, $(\alpha , \alpha ) =: m' = \min (L^*)$.
Dann liefert 
$$(\star _4) - (\star _2) : \ \ 
\sum _{x\in X} (x,\alpha )^2 (( x,\alpha )^2-1) = \frac{|X|mm'}{n} (mm'\frac{3}{n+2} -1 ) .$$
Da $(x,\alpha )\in \Z $ ist, ist die linke Seite eine nichtnegative Zahl,
also auch die rechte Seite, woraus die Behauptung folgt.
\hfill{q.e.d.}
\\

\noindent
\large
{\bf Ein Beispiel: Das Thompson-Smith Gitter.}
\normalsize
\\

\noindent
H\"aufig kann man die Kenntnis einer gro\ss en Untergruppe
$G\leq \Aut (L)$  der Automorphismengruppe
des Gitters $L$ benutzen, um zu zeigen, dass $L$ stark perfekt ist.
Die linke Seite von $(\star _4)$ ist n\"amlich ein $G$-invariantes homogenes
Polynom (in $\alpha $) vom Grad 4.
Hat $G$ keine anderen Invarianten von Grad 4 als das Quadrat der
invarianten quadratischen Form, dann ist die linke Seite von $(\star _4)$
ein Vielfaches von $(\alpha , \alpha )^2$ und damit nach Satz \ref{haupt} 
$L$ stark perfekt.
Nimmt man an, dass $-I_n\in G$ ist, so gilt sogar,
 dass alle $G$-Bahnen (und damit auch alle nicht leeren
Schichten $L_a$ von $L$) sph\"arische $5$-Designs sind.
Die Bedingung an die Invarianten vom Grad 4 von $G$ kann man leicht mit der 
Charaktertafel von $G$ nachpr\"ufen.
Sie ist z.B. f\"ur die 248-dimensionale Darstellung der sporadisch einfachen
Thompson Gruppe $Th$ erf\"ullt.
$\langle -I_n \rangle \times Th$ ist 
 Automorphismengruppe
eines (eindeutig bestimmten) geraden unimodularen Gitters,
dem {\sf Thompson Smith Gitter} $\Lambda _{248}$, der Dimension 248.
Dieses Gitter $\Lambda _{248}$ ist demnach stark perfekt also ein
lokales Maximum der Dichtefunktion.
Aus Satz \ref{min} folgert man, dass 
$$\min (\Lambda _{248}) \geq \sqrt{\frac{248+2}{3}} > 9 $$ 
also $\min (\Lambda _{248}) \geq 10$ ist. 
Diese Informationen erh\"alt man durch reine Charakterrechnung, ohne die 
Darstellung der Thompson Gruppe oder das Gitter $\Lambda _{248}$ 
explizit zu konstruieren.
Durch explizite Konstruktion der 248-dimensionalen Darstellung findet man
einen Vektor der L\"ange 12 in $\Lambda _{248}$ (in einem 1-dimensionalen 
invarianten Teilraum  der maximalen Untergruppe $(C_3\times G_2(3)):2$).
Jedoch ist immer noch offen, ob $\min (\Lambda _{248})$ nun 10 oder 12 ist.

\section{Theta-Reihen mit harmonischen Koeffizienten}
{\label{HARM}}

Mit Hilfe der Theorie der Modulformen kann man zeigen, dass gewisse extremale
Gitter stark perfekt und damit extrem sind.
Insbesondere sind alle extremalen geraden unimodularen Gitter in
den Sprungdimensionen $24l$ und in den Dimensionen $24l+8$ stark perfekt.
In Dimension 32 gibt es also mehr als 10 Millionen stark perfekte 
gerade unimodulare Gitter:

F\"ur ein harmonisches Polynom $P\in \Harm _t(2k) $ und ein gerades
Gitter $L$ von Stufe $N$ 
der Dimension $2k$ ist die {\sf Theta-Reihe von $L$ mit harmonischen Koeffizeinten $P$}
$$\theta _{L,P} := \sum _{x\in L} P(x) q^{(x,x)} $$
eine Modulform vom Gewicht $k+t$ 
zur Gruppe $\Gamma _0(N)$ mit dem Charakter $\chi _{k}$
(siehe z.B. \cite{Ebeling}, \cite{Miyake}).
Ist $t>0$, so ist $\theta _{L,P}$
 sogar eine Spitzenform.
Die Operation von $\Gamma _*(N) $ auf den Isometrieklassen von Gittern im
Geschlecht von $L$ \"ubersetzt sich f\"ur Theta-Reihen wie folgt:
Ist $N$ eine Primzahl, so ist 
$$ \theta _{L,P} + \theta _{L_{(N)},P} \in {\cal M}_{k+t} (\Gamma _*(N),\chi _{N,k}) $$
und
$$ \theta _{L,P} - \theta _{L_{(N)},P} \in {\cal S}_{k+t} (\Gamma _*(N),\chi _{N,k+2}) 
$$
(siehe \cite[Theorem 2.1]{BaV}).
F\"ur zusammengesetztes $N$ muss man $\pm 1$-Linear\-kombina\-tionen aller
partiellen dualen Gitter betrachten.
Durch Studium der entsprechenden Modulformen f\"ur $\Gamma _*(N)$
l\"a\ss t sich
zeigen, dass f\"ur kleine Grade $t\geq 1$ und extremale stark $N$-modulare Gitter $L$
beide Summen 
$$ \theta _{L,P} \pm \theta _{L_{(N)},P} = 0 $$ sind, also auch $\theta _{L,P}=0$.
Somit bilden alle Schichten 
$L_a$ von $L$ sph\"arische $t$-Designs.
Genauer findet man

\begin{Satz} (\cite[Corollary 3.1]{BaV})
Sei $L$ ein extremales $N$-modulares Gitter der Dimension $2k$.
Dann bilden die Schichten von $L$ sph\"arische $t$-Designs, gem\"a\ss  \ 
der folgenden Tabelle:
\begin{center}
\begin{tabular}{|c|c|c|c|c|c|c|}
\hline
$N$ & $1$ & $1$ & $2$ & $2$ & $3$ & $3$ \\
\hline
$k$ & $0 \mbox{ mod } 12$ & $4 \mbox{ mod } 12$ & $0 \mbox{ mod } 8$ & $2  \mbox{ mod } 8$   & $ 0  \mbox{ mod } 6$   &  $1  \mbox{ mod } 6$   \\
\hline
$t$ & $11$ & $7$ & $7$ & $5$ & $5$ & $5$ \\
\hline
\end{tabular}
\end{center}
Insbesondere ist $L$ in all diesen F\"allen  ein stark perfektes Gitter.
\end{Satz}

\underline{Beweis:} (f\"ur $N=1$).
Sei $L$ ein extremales gerades unimodulares Gitter der
Dimension $2k$ und $P \in \Harm _t(2k)$ ein harmonisches Polynom
vom Grad $t\geq 1$.
Dann ist $\min (L) = 2(1+\lfloor \frac{k}{12} \rfloor ) =:2m$ und 
$$
\theta _{L,P} = \sum _{j=m}^{\infty } (\sum _{x\in L_{2j}} P(x)) q^{2j} 
\in {\cal S}_{k+t} (\SL_2(\Z ) ).$$
Das Ideal der Spitzenformen
$$\bigoplus _{k=0}^{\infty  } {\cal S}_{k} (\SL_2(\Z )) = \Delta  \C [E_4 , E_6 ] $$
 f\"ur $\SL_2(\Z )$ ist ein Hauptideal in dem Polynomring
in $E_4 = \theta _{E_8}$ (vom Gewicht 4)
 und $E_6$ (vom Gewicht 6) erzeugt von der Spitzenform $\Delta = \Delta _1$ (vom Gewicht 12).
Da $\Delta $ mit $q^2$ beginnt, 
ist die mit $q^{2m}$ beginnende homogene Modulform
$\theta _{L,P} $ durch $\Delta ^m$ teilbar.
Das Gewicht von $\Delta ^m$ ist aber $12m = 12 + 12 \lfloor \frac{k}{12} \rfloor $.
Ist also $k \equiv 0 \pmod{12}$, so ist $\theta _{L,P} = 0$,
falls $t\leq 11$ ist und f\"ur 
$k \equiv 4 \pmod{12}$ ist $\theta _{L,P} = 0$, falls $t \leq 7$ ist.
Damit sind alle Schichten von $L$ sph\"arische 11- bzw. 7-Designs.
\hfill{q.e.d.}

Insbesondere liefern also die k\"urzesten Vektoren des Leech Gitters
ein 11-Design der Kardinalit\"at 
196560 auf $S^{23}$. Dies ist das einzige symmetrische 
11-Design mit $\leq 196560$ Punkten in Dimension 24
(\cite[Th\'eor\`eme 14.2]{Venkov}).

\section{Ungerade Gitter und deren Schatten}

Bisher haben wir nur  gerade Gitter $L$ betrachtet.
Extremalit\"at kann man auch f\"ur ungerade Gitter $L$ definieren,
wobei man beachten muss, dass $L$ nicht nur alle 
partiellen dualen Gitter festlegt, sondern auch sein
{\sf gerades Teilgitter}
$$L_{g} := \{ x\in L \mid (x,x) \in 2\Z \} $$
und dessen duales Gitter. 
$L_g$ ist der Kern der linearen Abbildung $L \to \F_2, x\mapsto (x,x) + 2\Z $, hat also Index 2 in $L$, falls $L$ ein ungerades Gitter ist.
Die Theta-Reihe des geraden Teilgitters ist
$$\theta _{L_g} (z) = \frac{1}{2} ( \theta _L(z) + \theta _L (z+1) ).$$

\begin{Definition}
Der {\sf Schatten} eines ungeraden Gitters $L$ ist 
$$S(L) := L_g^* - L^* .  $$
Ist $L$ gerade, so setzt man $S(L) = L^*$.
\end{Definition}

Der Schatten eines Gitters $L$ ist eine 
Restklasse nach $L^*$.
Er besteht aus den Vektoren $\frac{v}{2}$, wo
$v$ die {\sf charakteristischen Vektoren}
von $L$ durchl\"auft, also $(v,x) \equiv (x,x) \pmod{2}$ f\"ur alle $x\in L$. 
Mit der Theta-Transformationsformel ergibt sich die
 Theta-Reihe von $S(L)$ als
$$\theta _{S(L)} (z) = \sqrt{\det(L)} (\frac{i}{z} )^k \theta _L(1-\frac{1}{z}) .$$
Die Theta-Reihe eines ungeraden Gitters ist eine Modulform
f\"ur eine kleinere Gruppe als f\"ur die entsprechenden geraden Gitter.
Daher liefert die Modulformenbedingung f\"ur die Theta-Reihe 
alleine weniger scharfe Schranken als f\"ur gerade Gitter.
Als Ausgleich kann man aber auch noch die Theta-Reihe des
Schattens  betrachten, 
mit deren Hilfe man f\"ur dieselben Stufen
$N\in {\cal A} $ Extremalit\"at definieren kann.
Damit erh\"alt man dieselben  Schranken wie f\"ur
gerade Gitter.
Genauer gilt

\begin{Satz} (\cite[Theorem 1,2]{RaS}) {\label{ungeradeext}}
Sei $N\in {\cal A}$ und $$C_N := \perp _{d\mid N} \sqrt{d} \Z $$
ein stark $N$-modulares Gitter der Dimension $\sigma_0(N)$.
Ist $L$ ein stark $N$-modulares Gitter der Dimension $2k$,
welches rational \"aquivalent ist zu einer orthogonalen
Summe von Kopien von $C_N$, so gilt 
$$\min (L) \leq 2\lfloor \frac{k}{k_N} \rfloor + 2 ,$$
mit der Ausnahme $k = k_N - \frac{1}{2} \sigma_0(N) $, wo die Schranke
$3$ ist.
\end{Satz}

In Dimension $2k_N - \sigma_0(N)$ gibt es  genau ein stark $N$-modulares
Gitter $S^{(N)}$ von Minimum 3, die ``shorter version'' des entsprechenden
extremalen Gitters $E^{(N)}$ der Dimension $2k_N$.
\\

\noindent
\large
{\bf Gitter mit langem Schatten.}
\normalsize
\\

\nopagebreak
\noindent
Anstatt modulare Gitter mit gro\ss em Minimum zu suchen, kann
man auch nach solchen (ungeraden) Gittern $L$ fragen, f\"ur die 
die minimale L\"ange eines charakteristischen Vektors
(also die minimale L\"ange eines Vektors in $S(L)$)
m\"oglichst gro\ss \ ist.
Dies hat zwar f\"ur die Anwendung in der Informations\"ubertragung
keine Bedeutung, da der minimale Abstand von zwei Vektoren des
Schattens gleich $\min (L^*)$ ist, ist aber theoretisch
von Interesse.
Initiiert wurde diese Frage von N.D. Elkies (\cite{Elkies}, \cite{Elkies2}), 
der gezeigt hat,
dass das Standardgitter 
$\Z^n$, das einzige unimodulare Gitter der Dimension $n$ 
mit maximalem Schatten   
( $\min \{ (v,v) \mid v \in S(\Z ^n)\} = n/4$) ist.
Mit Hilfe der in \cite{RaS} angegebenen Formeln l\"a\ss t sich 
dies leicht auf die stark $N$-modularen Gitter in Satz \ref{ungeradeext}
verallgemeinern, wobei hier $C_N$ dieselbe Rolle spielt wie 
$\Z $ im Fall von $N =1$.
Ist $L$ n\"amlich ein stark $N$-modulares Gitter mit $\min (L) = 1 $,
so hat $L$ einen orthogonalen Summanden $C_N$.
Der Einfachheit halber nehmen wir im folgenden $N$ als ungerade an.

\begin{Satz}(\cite{Elkies}, \cite{Elkies2} f\"ur $N=1$,
\cite{preprint})\label{preprint}
Sei $N \in \{ 1,3,5,7,11,15 ,23 \}$ und 
$L$ ein stark $N$-modulares Gitter welches rational
\"aquivalent ist zu $C_N^l$. 
Dann gilt 
$$\mbox{$\min _0$} (S(L)) := 
\min \{ (v,v) \mid v\in S(L) \} = \frac{l \sigma_1(N)-8m}{4N} =: M^{(N)}(l,m)$$
f\"ur ein $m\in \Z _{\geq 0}$
\\
Ist $m=0$, so ist $L = C_N^l$.
\\
Ist $m=1$, so ist $L=C_N^{l_1} \perp M$ mit einem stark $N$-modularen Gitter $M$ von
Minimum $>1$ und mit $\min _0 (S(M)) = M^{(N)}(l-l_1,1)$.
Die Dimension von $M$ ist $\leq 2k_N - \dim (C_N )$.
\end{Satz}

\underline{Beweis:}
Aus den Formeln f\"ur die Theta-Reihe des Schattens 
in \cite[Corollary 3]{RaS} findet man die Formel
$M^{(N)}(l,m)$
f\"ur das m\"ogliche Minimum der Schattenvektoren.
 F\"ur $m=0$ ist die Theta-Reihe von $L$
eindeutig bestimmt und gleich der von $C_N^l$.
Die Vektoren der L\"ange 1 in $L$ erzeugen ein unimodulares 
Teilgitter $\cong \Z^l$, welches als orthogonaler Summand von
$L$ abspaltet 
$$L = \Z ^l \perp M .$$
Durch partielles Dualisieren findet man, dass $L$ f\"ur alle exakten
Teiler $m$ von $N$  einen orthogonalen Summanden $\sqrt{m} \Z^l$ hat.
Also ist $L \cong C_N^l$.
\\
F\"ur $m=1$ kann man nach Abspalten von orthogonalen Summanden $C_N$ annehmen, dass 
$\min (L) \geq 2$ ist. Dann ist wieder die Theta-Reihe von $L$ eindeutig bestimmt.
Die Anzahl der Vektoren der Norm 2 in $L$ ist 
$2l(\frac{24}{\sigma _1(N)} -l -1 )$ wodurch man die Schranke an $l$ erh\"alt.
Ist $l= \frac{24}{\sigma _1(N)} -1$, so ist $\min (L) = 3$ und $L= S^{(N)}$, 
das ``Shorter'' Gitter aus \cite[Table 1]{RaS}.
\hfill{q.e.d.}

Es ist ein offenes Problem, ob f\"ur jedes $m$, die Dimension eines 
stark $N$-modularen Gitters $L$ wie in Satz \ref{preprint} mit $\min (L) \geq 2$ und $\min _0(S(L)) \geq M^{(N)}(\dim(L)/\dim(C_N),m) $ nach oben
beschr\"ankt werden kann. 
F\"ur $m=0,1$ liefert Satz \ref{preprint} scharfe Schranken.
F\"ur $N=1$ hat Gaulter \cite{Gaulter} solche (nicht scharfen) Schranken f\"ur 
$m=2$ und $3$ angegeben.





\begin{thebibliography}{99}
\bibitem{Bachochurw} {\sl C. Bachoc:} 
Voisinage au sens de Kneser pour les r\'eseaux quaternioniens.
Comment. Math. Helvetici {\bf 70} (1995), 350-374.
\bibitem{Bachoc} {\sl C. Bachoc:} 
Applications of coding theory to the construction of modular lattices.
J. Comb. Th. Ser. A {\bf 78} (1997), 92-119.
\bibitem{BaN} {\sl C. Bachoc, G. Nebe:}
Extremal lattices of minimum 8 related to the Mathieu group $M_{22}$.
J. reine angew. Math. {\bf 494} (1998), 155-171.
\bibitem{BaV} {\sl C. Bachoc, B. Venkov:}
Modular forms, lattices and spherical designs.
In \cite{Mart}, 
87-112.
\bibitem{BaW} {\sl E.S. Barnes, G.E. Wall:} 
Some extreme forms defined in terms of abelian groups. 
J. Austral. Math. Soc. {\bf 1} (1959), 47-63.
\bibitem{SPLAG} {\sl J. H. Conway, N. J. A. Sloane:}
Sphere packings, lattices and groups. Springer, 3. Auflage (1998). 
\bibitem{Ebeling} {\sl W. Ebeling:}
Lattices and codes. Vieweg (1994).
\bibitem{Eichler} {\sl M. Eichler:}
Quadratische Formen und orthogonale Gruppen. Springer Grundlehren {\bf 63}, 
2. Auflage (1974)
\bibitem{Elkies} {\sl N.D. Elkies:}
A characterization of the $\Z^ n$ lattice.
Math. Res. Lett. 2 (1995), no. 3, 321-326.
\bibitem{Elkies2} {\sl N.D. Elkies:}
Lattices and codes with long shadows.
Math. Res. Lett. 2 (1995), no. 5, 643-651
\bibitem{Gaulter} {\sl M. Gaulter:}
Lattices without short characteristic vectors. Math. Res. Lett. 5 (1998), no. 3, 353-362.
\bibitem{King} {\sl O. King:} 
A mass formula for unimodular lattices with no roots.
    Mathematics of Computation, (to appear)
\bibitem{KKU} {\sl A.I. Kostrikin, I.A. Kostrikin, V.A. Ufnarovskii:}
Invariant lattices of type $G_2$ and their automorphism groups. 
Proc. Steklov Inst. Math. {\bf 3} (1985), 85-105.
\bibitem{Mart} {\sl J. Martinet:} (Herausgeber)
R\'eseaux euclidiens, designs sph\'eriques et formes modulaires.
L'Ens. Math. Monographie {\bf 37} (2001).
\bibitem{MOS} {\sl C. L. Mallows, A. M. Odlysko, N. J. A. Sloane:}
Upper bounds for modular forms, lattices and codes.
J. Alg. {\bf 36} (1975), 68-76.
\bibitem{Miyake} {\sl T. Miyake:} Modular Forms. Springer (1989).
\bibitem{dim24} {\sl G. Nebe:}
Endliche rationale Matrixgruppen vom Grad 24.
Dissertation, 
Verlag der Augustinus Buchhandlung,
Aachener Beitr\"age zur Mathematik {\bf 12} (1995).
Englische Kurzfassung:
Finite subgroups of $GL_{24}(\Q )$. 
Exp. Math. {\bf 5} (3) (1996), 163-195.
\bibitem{dim2531} {\sl G. Nebe:}
Finite subgroups of $GL_n(\Q )$ for $25\leq n \leq 31$.
Comm. Alg. {\bf 24} (7), (1996), 2341-2397.
\bibitem{quat} {\sl G. Nebe:}
Finite quaternionic matrix groups. 
Representation Theory {\bf 2}, (1998), 106-223.
\bibitem{cyclo} {\sl G. Nebe:}
Some cyclo-quaternionic lattices. J. Alg. {\bf 199} (1998), 472-498.
\bibitem{seoul} {\sl G. Nebe:}
Construction and investigation of lattices with matrix groups.
in Integral Quadratic Forms and Lattices,
herausgegeben von 
 Myung-Hwan Kim, John S. Hsia, Y. Kitaoka, R. Schulze-Pillot, 
Contemporary Mathematics {\bf 249} (1999), 205-220.
\bibitem{preprint} {\sl G. Nebe:}
Strongly modular lattices with long shadow. (Preprint 2002)
\bibitem{NeP} {\sl G. Nebe, W. Plesken:} Finite rational matrix groups.
AMS Memoirs {\bf 116} (556) (1995).
\bibitem{database} {\sl G. Nebe, N.J.A. Sloane:} 
A database of lattices. Internetsammlung unter
 www.research.att.com/$\sim$njas/lattices
\bibitem{NeV} {\sl G. Nebe, B. Venkov:} 
 Non-existence of extremal lattices in certain genera of modular lattices.
J. Number Theory, {\bf 60} (2) (1996), 310-317. 
\bibitem{Siegel} {\sl G. Nebe, B. Venkov:} 
On Siegel modular forms of weight 12.
J. reine und angew. Mathematik {\bf 531} (2001), 49-60.
\bibitem{queb32q} {\sl H.-G. Quebbemann:} 
A construction of integral lattices. Mathematika {\bf 31} (1984), 137-140.
\bibitem{queb32} {\sl H.-G. Quebbemann:} 
 Lattices with Theta Functions for $G(\sqrt{2})$ and Linear Codes.
J. Algebra {\bf 105} (1987),  443-450.
\bibitem{queb1} {\sl H.-G. Quebbemann:} 
Modular lattices in euclidean spaces. J. Number Th. {\bf 54} (1995), 190-202.
\bibitem{queb2} {\sl H.-G. Quebbemann:} 
Atkin-Lehner eigenforms and strongly modular lattices.
L'Ens. Math. {\bf 43} (1997), 55-65.
\bibitem{PohstZas} {\sl M. Pohst, H. Zassenhaus:}
Algorithmic algebraic number theory. Cambridge University Press (1989).
\bibitem{RaS} {\sl E.M. Rains, N.J.A. Sloane:}
The shadow theory of modular and unimodular lattices. 
J. Number Th. {\bf 73} (1998), 359-389.
\bibitem{SchaHe} {\sl R. Scharlau, B. Hemkemeier:}
Classification of integral lattices with large class number.
Math. Comp. {\bf 67} (1998), 737-749.
\bibitem{SchaSchuPi} {\sl R. Scharlau, R. Schulze-Pillot:}
Extremal lattices. In 
Algorithmic algebra and number theory. 
Herausgegeben von B. H. Matzat, G. M. Greuel, G. Hiss. 
Springer (1999), 139-170.
Preprint erh\"altlich unter 
www.matha.mathematik.uni-dortmund.de/preprints/welcome.html
\bibitem{SchV} {\sl R. Scharlau, B. Venkov:}
The genus of the Barnes-Wall lattice.
Comm. Math. Helvetici {\bf 69} (1994), 322-333.
\bibitem{SchV2} {\sl R. Scharlau, B. Venkov:}
The genus of the Coxeter-Todd lattice.
Preprint 95-7
erh\"altlich unter 
www.matha.mathematik.uni-dortmund.de/preprints/welcome.html
\bibitem{Serre} {\sl J. P. Serre:}
A course in arithmetic. Springer GTM {\bf 7} (1973).
\bibitem{ShT} {\sl G.C. Shephard, J.A. Todd:} 
Finite unitary reflection groups. Canad. J. Math. {\bf 6} (1954), 274-304.
\bibitem{CLSiegel} {\sl C.L. Siegel:} Berechnung von Zetafunktionen an ganzzahligen Stellen. Nachr. Akad. Wiss. G\"ottingen {\bf 10} (1969), 87-102.
\bibitem{Tits} {\sl J. Tits:}
Quaternions over $\Q[\sqrt{5}]$, Leech's lattice and the sporadic group of Hall-Janko.
J. Alg. {\bf 63} (1980), 56-75.
\bibitem{Venkov} {\sl B. Venkov:}
R\'eseaux et designs sph\'eriques. 
In \cite{Mart} 10-86.
\bibitem{Voronoi} {\sl G. Voronoi:}
Nouvelles applications des param\`etres continues \`a la th\'eorie des formes
quadratiques: 1. Sur quelques propriet\'es des formes quadratiques parfaites.
J. Reine Angew. Math. {\bf 133} (1908) 97-178.
\end{thebibliography}
\end{document}